\documentclass[10pt]{amsart}
\usepackage{latexsym, amsmath}

\renewcommand{\epsilon}{\varepsilon}
\newcommand{\newsection}[1]
{\subsection{#1}\setcounter{theorem}{0} \setcounter{equation}{0}
\par\noindent}

\newtheorem{theorem}{Theorem}

\newtheorem{lemma}[theorem]{Lemma}

\newtheorem{proposition}[theorem]{Proposition}

\newcommand{\cd}{\, \cdot\, }

\newcommand{\Rplus}{{\Bbb R}_+}
\newcommand{\R}{{\mathbb R}}
\newcommand{\ext}{{\R^n\backslash\mathcal{K}}}
\newcommand{\extiv}{{\R^4\backslash\mathcal{K}}}

\begin{document}

\title[Global existence for Dirichlet-wave equations]
{Global existence for Dirichlet-wave equations with quadratic
nonlinearities in high dimensions}
\thanks{The authors were supported in part by the NSF}

\author{Jason Metcalfe}
\address{School of Mathematics, Georgia Institute of Technology,
Atlanta, GA  30332-0160}

\author{Christopher D. Sogge}
\address{Department of Mathematics, Johns Hopkins University,
Baltimore, MD 21218}

\maketitle

\newsection{Introduction}

In this paper, we provide a proof of global existence of solutions to
quasilinear wave equations with quadratic nonlinearities exterior to
nontrapping obstacles.  Specifically,  let $\mathcal{K}$ be a compact,
nontrapping obstacle with smooth boundary.  We will then be looking
for solutions to
\begin{equation}\label{1.1}
\begin{cases}
\Box u = Q(du,d^2u),\quad (t,x)\in \R_+\times\ext\\
u(0,\cd)=f,\quad \partial_t u(0,\cd)=g\\
u(t,\cd)|_{\partial\mathcal{K}}=0
\end{cases}
\end{equation}
exterior to $\mathcal{K}$.  Here
$$\Box=(\Box_{c_1},\Box_{c_2},\dots,\Box_{c_D})$$
is a vector-valued d'Alembertian with
$$\Box_{c_I}=\partial^2_t-c_I^2\Delta$$
and $c_I>0$ for $I=1,2,\dots,D$.

Letting $\partial_0=\partial_t$ when convenient, we can expand our
quadratic, quasilinear forcing term $Q$ as follows
\begin{equation}\label{1.2}
Q^I(du,d^2u)=\sum_{\substack{0\le j,k\le
n\\ 1\le J,K\le D}}A^{I,jk}_{JK}\partial_j u^J \partial_k u^K 
+ \sum_{\substack{0\le j,k,l\le n\\ 1\le J,K\le D}}B^{IJ,jk}_{K,l}
\partial_l u^K \partial_j\partial_k u^J,\quad 1\le I\le D.
\end{equation}
In order that we might apply the local existence results of
Keel-Smith-Sogge \cite{KSS} and in order to help guarantee
hyperbolicity, we assume the following symmetry condition
\begin{equation}\label{1.3}
B^{IJ,jk}_{K,l}=B^{JI,jk}_{K,l}=B^{IJ,kj}_{K,l}.
\end{equation}

To solve \eqref{1.1}, one must assume that the Cauchy data $(f,g)$
satisfy certain compatibility conditions.  Such conditions are
well-known, and for further detail, we refer the reader to
\cite{KSS}.  Briefly, if we let $J_ku=\{\partial_x^\alpha u\,:\, 0\le
|\alpha|\le k\}$, we can write $\partial_t^k
u(0,\cd)=\psi_k(J_kf,J_{k-1}g)$, $0\le k\le m$, where $u$ is any formal
$H^m$ solution to \eqref{1.1} and $m$ is fixed.  The $\psi_k$ are
called compatibility functions and depend on $Q$, $J_kf$, and
$J_{k-1}g$.  The compatibility condition for \eqref{1.1} with
$(f,g)\in H^m\times H^{m-1}$ requires that $\psi_k$ vanish on
$\partial\mathcal{K}$ when $0\le k\le m-1$.  Additionally, we say that
$(f,g)\in C^\infty$ satisfy the compatibility condition to infinite
order if the above condition holds for all $m$.

The main result of this paper asserts that such systems of multiple
speed, Dirichlet-wave equations admit global solutions.

\begin{theorem}\label{theorem1.1}  Assume $n\ge 4$, and
let $\mathcal{K}$, $Q$ and $\Box$ be as above.  Suppose further that
$(f,g)\in C^\infty(\ext)$ satisfy the compatibility conditions to
infinite order.  Then, there is a constant $\varepsilon_0>0$ and an
integer $N>0$ so that for all $\varepsilon \le \varepsilon_0$, if
\begin{equation}\label{1.4}
\sum_{|\alpha|\le N} \|\langle x\rangle^{|\alpha|}\partial_x^\alpha
f\|_2 + \sum_{|\alpha|\le N-1}\|\langle
x\rangle^{1+|\alpha|}\partial^\alpha_x g\|_2\le \varepsilon,
\end{equation}
then \eqref{1.1} has a unique global solution $u\in C^\infty([0,\infty)\times\ext)$.
\end{theorem}

Additionally, we note that the proof of the theorem would allow any
forcing term $F(du,d^2u)$ vanishing to second order and linear in
$d^2u$.

Global existence of solutions to boundaryless wave equations of the
form \eqref{1.1} was first shown by H\"ormander and Klainerman (see,
e.g., \cite{S}).  A recent paper of Hidano \cite{Hidano}
explores an alternate method of proof that admits the multiple speed
setting.

In the obstacle setting, \eqref{1.1} was first considered by 
Shibata-Tsutsumi \cite{tsutsumi} and was shown to have global existence in
spatial dimensions $n\ge 6$.  Hayashi \cite{Hayashi} was able to prove
global existence exterior to a ball in all spatial dimensions $n\ge
4$.  A result similar to Theorem \ref{theorem1.1} was shown by the
first author \cite{Metcalfe} for semilinear equations.

In the case of $n=3$, solutions to \eqref{1.1} exterior to certain
obstacles were studied by Keel-Smith-Sogge \cite{KSS2,KSS3}, the
authors \cite{MS}, and Metcalfe-Nakamura-Sogge \cite{MNS}.  As in
these works, we will be using the exterior domain analog of
Klainerman's method of commuting vector fields \cite{knull} as
developed by Keel-Smith-Sogge \cite{KSS3}.  In particular, we restrict
our attention to the invariant vector fields that are admissible for
the obstacle setting, $\{L,Z\}$, where $Z$ represents the generators
of the space-time translations and spatial rotations
$$Z=\{\partial_i, x_j\partial_k - x_k\partial_j\}, \quad 0\le i\le
n,\quad 1\le j,k\le n$$
and where $L$ is the scaling vector field
$$L=t\partial_t + r\partial_r.$$
Here and in what follows, $r=|x|$.  We also set
$$\Omega=\{x_j\partial_k-x_k\partial_j\},\quad 1\le j,k\le n$$
to be the set of generators of spatial rotations.

The main new approach in this paper versus \cite{MS} is the techniques
used to handle the boundary terms that necessarily arise when studying
obstacle problems.  In \cite{MS}, these were handled using Huygens'
principle.  In the current setting, we develop simple local bounds for
solutions to the Minkowski wave equation using the fundamental
solution.  We then use local energy decay and techniques of
Smith-Sogge \cite{SS} to reduce to this case.

Also, as in \cite{KSS3}, we will be using a class of weighted
$L^2_tL^2_x$ estimates where the weight is a negative power of
$\langle x\rangle=\langle r\rangle=\sqrt{1+r^2}$.  Such estimates
allow one to take advantage of the decay in $|x|$ which is much easier
to prove in the obstacle setting than the more traditional decay in
$t$.  These estimates were first developed for even spatial dimensions
by the first author in \cite{Metcalfe}.  The proof relied on a local
Minkowski version developed by Smith-Sogge \cite{SS} and on other weighted
estimates in \cite{Metcalfe2}.  Local versions of these weighted
$L^2_tL^2_x$ estimates were originated in the obstacle setting 
using different techniques by
Burq \cite{B}.  Burq's estimates relied on rather weak hypotheses,
namely the existence of certain resolvant bounds.  Since these
resolvant bounds are implied by the local energy decay that we discuss
next, we will assume Burq's bounds when convenient.

By a simple scaling argument, we may and will assume throughout that
$$\mathcal{K}\subset \{|x|\le 1\}.$$
The nontrapping assumption on the geometry of the obstacle, which states that
there is a $T_R$ such that no geodesic of length $T_R$ is completely
contained in $\{|x|\le R\}\cap \ext,$
allows us to refer to well-known local energy decay
estimates.  In particular, if $u$ is a solution to the homogeneous
wave equation
\begin{equation}\label{1.5}
\begin{cases}
(\partial_t^2-\Delta)u(t,x)=0\\
u(0,\cd)=f,\quad \partial_tu(0,\cd)=g\\
u(t,x)=0, \quad x\in\partial\mathcal{K}
\end{cases}
\end{equation}
and if the Cauchy data $f,g$ are assumed to vanish for $|x|>4$, then
there is a constant $c>0$ so that
\begin{equation}\label{1.6}
\Bigl(\int_{\{x\in\ext : |x|<4\}} |u'(t,x)|^2\:dx\Bigr)^{1/2}\le C
e^{-ct}\Bigl(\|\nabla_x f\|_2 + \|g\|_2\Bigr) 
\end{equation}
if $n$ is odd.  Here $'=(\partial_t,\nabla_x)$ is the space-time
gradient.  We refer
the reader to Taylor \cite{Taylor}, Lax-Phillips \cite{LP}, Vainberg
\cite{Vainberg}, Morawetz-Ralston-Strauss \cite{MRS}, Strauss
\cite{Strauss}, and Morawetz \cite{M}.

In even spatial dimensions $n$, we have the weaker decay
\begin{equation}\label{1.7}
\Bigl(\int_{\{x\in\ext : |x|<4\}} |u'(t,x)|^2\:dx\Bigr)^{1/2}\le C
t^{-(n-1)}\Bigl(\|\nabla_x f\|_2 + \|g\|_2\Bigr). 
\end{equation}
See Ralston \cite{Ralston2}.  We also refer the reader to Melrose
\cite{Mel} and Strauss \cite{Strauss}.  We will not require the
additional decay \eqref{1.6} and will only use \eqref{1.7} throughout.

One of the advantages of the proof that we shall use is
that the argument can easily be altered to allow for the necessary
loss of regularity in the right sides of \eqref{1.6} and \eqref{1.7}
if the exterior domain contains trapped rays.  The necessity of such a
loss was shown by Ralston \cite{ralston}, and in $n=3$, Ikawa
\cite{Ikawa1, Ikawa2} was able to show a version of \eqref{1.6} with a
loss of regularity for certain exterior domains that contain
hyperbolic trpaped rays.  In \cite{MS,MNS}, considerations were taken
to establish existence results in the presence of such geometries.

This paper is organized as follows.  In the next section, we collect
the $L^2$ energy estimates that we will require.  These are $n\ge 4$
analogs of those developed by the authors in \cite{MS}, and the proofs
of these results extend trivially to the more general setting.  In \S
3, we prove the necessary weighted $L^2_tL^2_x$ estimates.  As
mentioned previously, these follow easily from the estimates in
\cite{B} and \cite{Metcalfe} and are higher dimensional analogs of the
estimates of Keel-Smith-Sogge \cite{KSS2, KSS3}.  In \S 4, we state a
few Sobolev-type results.  These are exterior domain analogs of
results proven and used by Klainerman \cite{knull}, Klainerman-Sideris
\cite{KS}, Sideris \cite{Si}, Sideris-Tu \cite{Si3}, and
Hidano-Yokoyama \cite{HY, HY2}.  The extension of these estimates to
the exterior domain follows exactly as in \cite{MNS}.  In \S 5, we
provide proofs of the estimates for the boundary terms.  Finallly, in
\S 6, we set up a continuity argument and use these estimates to prove
global existence.

%%%%%%%%%%%%%%%%%%%%%%%%%%%%%%%%%%%%%%%%%%%%%%%%%%%%%%%%%%%%%%%%
%%%%%%%%%%%%%% SECTION 2 %%%%%%%%%%%%%%%%%%%%%%%%%%%%%%%%%%%%%%%
%%%%%%%%%%%%%%%%%%%%%%%%%%%%%%%%%%%%%%%%%%%%%%%%%%%%%%%%%%%%%%%%

\newsection{Energy Type Estimates}

In this section, we collect the energy estimates that we shall
require.  Unless stated otherwise, the proofs of these estimates can
be found in \cite{MS} for the $n=3$ case.  These arguments, however,
extend to general spatial dimensions $n\ge 2$ trivially.

Specifically, we will be concerned with solutions $u\in C^\infty
(\Rplus\times\ext)$ of the Dirichlet-wave equation
\begin{equation}\label{2.1}
\begin{cases}
\Box_\gamma u=F\\
u|_{\partial\mathcal{K}}=0\\
u|_{t=0}=f,\quad \partial_t u|_{t=0}=g
\end{cases}
\end{equation}
where
$$(\Box_\gamma u)^I=(\partial_t^2-c_I^2\Delta)u^I
+\sum_{J=1}^D\sum_{j,k=0}^n
\gamma^{IJ,jk}(t,x)\partial_j\partial_k u^J,\quad 1\le I\le D.$$
We shall assume that the $\gamma^{IJ,jk}$ satisfy the symmetry
conditions
\begin{equation}\label{2.2}
\gamma^{IJ,jk}=\gamma^{JI,jk}=\gamma^{IJ,kj}
\end{equation}
as well as the size condition
\begin{equation}\label{2.3}
\sum_{I,J=1}^D \sum_{j,k=0}^n \|\gamma^{IJ,jk}(t,x)\|_{\infty}\le
\delta
\end{equation}
for $\delta$ sufficiently small (depending on the wave speeds).
 The energy estimate will involve bounds for the gradient of the
perturbation terms
$$\|\gamma'(t,\cd)\|_\infty = \sum_{I,J=1}^D\sum_{j,k,l=0}^n
\|\partial_l \gamma^{IJ,jk}(t,\cd)\|_\infty,$$ and the energy form
associated with $\Box_\gamma$, $e_0(u)=\sum_{I=1}^D e_0^I(u)$,
where
\begin{multline}\label{2.4}
e_0^I(u)=(\partial_0 u^I)^2+\sum_{k=1}^n c_I^2(\partial_k u^I)^2
\\+2\sum_{J=1}^D\sum_{k=0}^n
\gamma^{IJ,0k}\partial_0u^I\partial_ku^J -
\sum_{J=1}^D\sum_{j,k=0}^n
\gamma^{IJ,jk}\partial_ju^I\partial_ku^J.
\end{multline}

The most basic estimate will lead to a bound for
$$E_M(t)=E_M(u)(t)=\int \sum_{j=0}^M e_0(\partial^j_t
u)(t,x)\:dx.$$

\begin{lemma}\label{lemma2.1}
Fix $M=0,1,2,\dots$, and assume that the perturbation terms
$\gamma^{IJ,jk}$ are as above.  Suppose also that $u\in C^\infty$
solves \eqref{2.1} and for every $t$, $u(t,x)=0$ for large $x$.
Then there is an absolute constant $C$ so that
\begin{equation}\label{2.5}
\partial_t E^{1/2}_M(t)\le C\sum_{j=0}^M \|\Box_\gamma
\partial_t^ju(t,\cd)\|_2+C\|\gamma'(t,\cd)\|_\infty E^{1/2}_M(t).
\end{equation}
\end{lemma}

Before stating the next result, let us introduce some notation. If
$P=P(t,x,D_t,D_x)$ is a differential operator, we shall let
$$[P,\gamma^{kl}\partial_k\partial_l]u=\sum_{1\le I,J\le
D}\sum_{0\le k,l\le n} |[P,\gamma^{IJ,kl}\partial_k\partial_l]
u^J|.$$

In order to allow the above energy estimate to include the
more general vector fields $L, Z$, we will need to use a variant
of the scaling vector field $L$.  We fix a bump function $\eta\in
C^\infty(\R^n)$ with $\eta(x)=0$ for $x\in \mathcal{K}$ and
$\eta(x)=1$ for $|x|>1$.  Then, set $\tilde{L}=\eta(x)r\partial_r
+ t\partial_t$. Using this variant of the scaling vector field and
an elliptic regularity argument, one can establish

\begin{proposition}\label{proposition2.2}
Suppose that the constant in \eqref{2.3} is small.  Suppose
further that
\begin{equation}\label{2.6}
\|\gamma'(t,\cd)\|_\infty \le \delta/(1+t),
\end{equation}
and
\begin{multline}\label{2.7}
\sum_{\substack{j+\mu\le
N_0+\nu_0\\\mu\le\nu_0}}\left(\|\tilde{L}^\mu\partial^j_t\Box_\gamma
u(t,\cd)\|_2+\|[\tilde{L}^\mu\partial_t^j,\gamma^{kl}\partial_k\partial_l]u(t,\cd)\|_2\right)\\
\le \frac{\delta}{1+t}\sum_{\substack{j+\mu\le N_0+\nu_0\\
\mu\le\nu_0}}\|\tilde{L}^\mu\partial_t^j
u'(t,\cd)\|_2+H_{\nu_0,N_0}(t),
\end{multline}
where $N_0$ and $\nu_0$ are fixed.  Then
\begin{multline}\label{2.8}
\sum_{\substack{|\alpha|+\mu\le N_0+\nu_0\\
\mu\le\nu_0}}\|L^\mu\partial^\alpha u'(t,\cd)\|_2 \\ \le
C\sum_{\substack{|\alpha|+\mu\le N_0+\nu_0-1\\
\mu\le\nu_0}}\|L^\mu\partial^\alpha\Box u(t,\cd)\|_2 +
C(1+t)^{A\delta}\sum_{\substack{\mu+j\le N_0+\nu_0\\ \mu\le
\nu_0}}\left(\int e_0(\tilde{L}^\mu\partial_t^j
u)(0,x)\:dx\right)^{1/2}\\
+C(1+t)^{A\delta}\Bigl(\int_0^t \sum_{\substack{|\alpha|+\mu\le
N_0+\nu_0-1\\\mu\le\nu_0-1}}\|L^\mu\partial^\alpha \Box
u(s,\cd)\|_2\:ds +\int_0^t H_{\nu_0,N_0}(s)\:ds\Bigr)\\
+C(1+t)^{A\delta}\int_0^t \sum_{\substack{|\alpha|+\mu\le
N_0+\nu_0\\ \mu\le \nu_0-1}}\|L^\mu\partial^\alpha
u'(s,\cd)\|_{L^2(|x|<1)}\:ds,
\end{multline}
where the constants $C$ and $A$ are absolute constants.
\end{proposition}

In practice $H_{\nu_0,N_0}(t)$ will involve weighted $L^2_x$ norms
of $|L^\mu \partial^\alpha u'|^2$ with $\mu+|\alpha|$ much smaller
than $N_0+\nu_0$, and so the integral involving $H_{\nu_0,N_0}$
can be dealt with using an inductive argument and the weighted
$L^2_tL^2_x$ estimates of the subsequent section.

In proving our existence results for \eqref{1.1}, a key step
will be to obtain a priori $L^2$-estimates involving $L^\mu
Z^\alpha u'$.  Begin by setting
\begin{equation}\label{2.9}
Y_{N_0,\nu_0}(t)=\int \sum_{\substack{|\alpha|+\mu\le
N_0+\nu_0\\\mu\le\nu_0}}e_0(L^\mu Z^\alpha u)(t,x)\:dx.
\end{equation}
We, then, have the following proposition which shows how the
$L^\mu Z^\alpha u'$ estimates can be obtained from the ones
involving $L^\mu\partial^\alpha u'$.
\begin{proposition}\label{proposition2.3}
Suppose that the constant $\delta$ in \eqref{2.3} is small and
that \eqref{2.6} holds.  Then,
\begin{multline}\label{2.10}
\partial_t Y_{N_0,\nu_0}\le C Y^{1/2}_{N_0,\nu_0} \sum_{\substack{
|\alpha|+\mu\le N_0+\nu_0\\ \mu\le\nu_0}} \|\Box_\gamma L^\mu Z^\alpha
u(t,\cd)\|_2 + C \|\gamma'(t,\cd)\|_\infty Y_{N_0,\nu_0} \\
+C \sum_{\substack{|\alpha|+\mu\le N_0+\nu_0+1\\ \mu\le
\nu_0}} \|L^\mu \partial^\alpha u'(s,\cd)\|^2_{L^2(|x|<1)}.
\end{multline}
\end{proposition}

%%%%%%%%%%%%%%%%%%%%%%%%%%%%%%%%%%%%%%%%%%%%%%%%%%%%%%%%%%%%%%%%
%%%%%%%%%%%%%% SECTION 3 %%%%%%%%%%%%%%%%%%%%%%%%%%%%%%%%%%%%%%%
%%%%%%%%%%%%%%%%%%%%%%%%%%%%%%%%%%%%%%%%%%%%%%%%%%%%%%%%%%%%%%%%

\newsection{$\mathbf{L^2_tL^2_x}$ Estimates}

As in Keel-Smith-Sogge \cite{KSS2, KSS3}, we will require a class of
weighted $L^2_tL^2_x$ estimates.  They will be used, for example, to
control the local $L^2$ norms such as the last term in \eqref{2.10}.
For convenience, allow $\Box=\partial_t^2-\Delta$ to denote the unit
speed, scalar d'Alembertian for the remainder of the section.  The
transition to the general case is straightforward.  Also, set
$$S_T=\{[0,T]\times\ext\}$$
to be the time strip of height $T$ in $\R_+\times\ext$.  Here we will
study solutions of the wave equation with vanishing initial data.  In
the sequel, we will reduce to this case.

We first note that if $u$ is a solution to 
\begin{equation}\label{3.1}
\begin{cases}
\Box u(t,x)=F(t,x)+G(t,x),\quad (t,x)\in \R_+\times\ext\\
u(t,x)|_{\partial\mathcal{K}}=0\\
u(t,x)=0,\quad t<0
\end{cases}
\end{equation}
where $G(t,x)=0$ for $|x|>2$,
then we have that
\begin{equation}\label{3.2}
\|u'\|_{L^2_tL^2_x([0,t]\times \{|x|<2\})} \le C \int_0^t
\|F(s,\cd)\|_2\:ds + C\|G\|_{L^2_tL^2_x([0,t]\times\ext)}.\end{equation}
Indeed, \eqref{3.2} was shown to follow from certain resolvant
estimates in Burq \cite{B} (Theorem 3).  Since the local energy decay estimates
\eqref{1.6} and \eqref{1.7} imply these resolvant estimates, we assume \eqref{3.2}.

Since $[\partial_t,\Box]=0$ and since $\partial_t$ preserves the
boundary condition, \eqref{3.2} holds with $u$ replaced by
$\partial_t^j u$ and $F, G$ replaced by $\partial_t^j F, \partial_t^j
G$ respectively for any $j=0,1,2,\dots$.  By elliptic regularity (see
Lemma 2.3 of \cite{MS}), it follows that
\begin{multline}\label{3.3}
\sum_{|\alpha|\le N} \|\partial^\alpha
u'\|_{L^2_tL^2_x([0,t]\times\{|x|<2\})} \le C\sum_{|\alpha|\le N}
\int_0^t \|\partial^\alpha F(s,\cd)\|_2\:ds 
\\+ C\sum_{|\alpha|\le N-1}
\|\partial^\alpha F\|_{L^2_tL^2_x([0,t]\times\ext)}
+C\sum_{|\alpha|\le N} \|\partial^\alpha G\|_{L^2_tL^2_x([0,t]\times\ext)}
\end{multline}
if $G$ is as above.  Moreover, using an inductive argument in $\nu_0$, we can prove
\begin{lemma}
\label{lemma_local}  Suppose $n\ge 3$, and suppose that $\mathcal{K}$ is nontrapping.
Let $u$ be a solution to \eqref{3.1}, and suppose $G(t,x)=0$ for $|x|>2$.  Then, for
any integers $\nu_0,N\ge 0$, we have
\begin{multline}
\label{3.local}
\sum_{\substack{|\alpha|+\mu\le N+\nu_0\\\mu\le \nu_0}} \|L^\mu \partial^\alpha u'\|_{L^2_tL^2_x([0,t]\times
\{|x|<2\})} \le C \sum_{\substack{|\alpha|+\mu\le N+\nu_0\\\mu\le \nu_0}} \int_0^t \|L^\mu \partial^\alpha F(s,\cd)\|_2\:ds
\\+C\sum_{\substack{|\alpha|+\mu\le N+\nu_0-1\\\mu\le \nu_0}} \|L^\mu \partial^\alpha F\|_{L^2_tL^2_x([0,t]\times 
\R^n\backslash\mathcal{K})} + C\sum_{\substack{|\alpha|+\mu\le N+\nu_0\\\mu\le\nu_0}} \|L^\mu \partial^\alpha G\|_{
L^2_tL^2_x([0,t]\times\R^n\backslash\mathcal{K})}.
\end{multline}
\end{lemma}

\noindent{\em Proof of Lemma \ref{lemma_local}:}  We will indeed use induction on $\nu_0$ where \eqref{3.3} serves
as the base case $\nu_0=0$.  We now proceed under the the assumption that \eqref{3.local} holds 
for any $N$ with $\nu_0$ replaced
by $\nu_0-1$.

Letting $\tilde{L}$ be as in the previous section, we see that the left side of \eqref{3.local} is dominated by
$$\sum_{|\alpha|\le N} \|L^{\nu_0-1} \partial^\alpha (\tilde{L}u)'\|_{L^2_tL^2_x([0,t]\times \{|x|<2\})}
+\sum_{\substack{|\alpha|+\mu\le N+\nu_0\\\mu\le \nu_0-1}} \|L^\mu \partial^\alpha u'\|_{L^2_tL^2_x([0,t]\times
\{|x|<2\})}.$$
By the inductive hypothesis, the second term is trivially controlled by the right side of \eqref{3.local}.

For the first term, we note that
\begin{multline}\label{commutator}
\Box(\tilde{L}u)=\tilde{L}\Box u + [\Box,\tilde{L}]u  = \tilde{L}\Box u + 2\Box u + [\Box, (1-\eta(x))r\partial_r]u
\\=\Bigl(\tilde{L}\Box u + 2\Box u\Bigr) + \Bigl(-(\Delta \eta)r\partial_r u
-2\partial_r \eta \partial_r u - 2r\nabla\eta\cd\nabla_x (\partial_r u) + 2(1-\eta) \partial_r^2 u\Bigr).
\end{multline}
Notice, in particular, that the second grouping of terms are all supported in $|x|<1$.  Thus, if we apply
the inductive hypothesis to $\tilde{L}u$, it follows that the left side of \eqref{3.local} is bounded by
the right side of \eqref{3.local} plus
$$\sum_{\substack{|\alpha|+\mu\le N+\nu_0\\\mu\le \nu_0-1}} \|L^\mu \partial^\alpha u'\|_{L^2_tL^2_x([0,t]
\times \{|x|<1\})}.$$
This last term, using the inductive hypothesis, is also easily seen to be controlled by the right side
of \eqref{3.local} which completes the proof.
\qed

We will also require the associated global results in Minkowski
space.  Let $v$ be a solution to the boundaryless wave equation
\begin{equation}\label{3.4}
\begin{cases}
\Box v(t,x)=F(t,x)+G(t,x),\quad (t,x)\in \R_+\times \R^n\\
v(t,x)=0,\quad t<0.
\end{cases}
\end{equation}
We, then, have the following result of the first author
\cite{Metcalfe} (Proposition 2.2 and Proposition 2.7)
\begin{lemma}\label{lemma3.1}
Suppose $n\ge 4$, and let $v$ be a solution to \eqref{3.4}.  If
$G(s,x)=0$ for $|x|>2$, then
\begin{multline}\label{3.5}
\sum_{\substack{|\alpha|+\mu\le N+\nu\\\mu\le \nu}} \|\langle
r\rangle^{-(n-1)/4} L^\mu Z^\alpha v'\|_{L^2_tL^2_x([0,t]\times\R^n)}
\le C\int_0^t \sum_{\substack{|\alpha|+\mu\le N+\nu\\\mu\le \nu}}
\|L^\mu Z^\alpha F(s,\cd)\|_2\:ds
\\+C\sum_{\substack{|\alpha|+\mu\le N+\nu\\\mu\le \nu}} \|L^\mu
\partial^\alpha G\|_{L^2_tL^2_x([0,t]\times\R^n)}
\end{multline}
for any $N,\nu\ge 0$.
\end{lemma}

A global estimate for the Dirichlet-wave equation 
will follow from the local estimate \eqref{3.local} and the global
Minkowski estimate \eqref{3.5}.  In particular, we have the following
$n\ge 4$ analog of Theorem 6.3 of Keel-Smith-Sogge \cite{KSS3}.
\begin{proposition}\label{prop3.2}
Fix $N_0$ and $\nu_0$.  Suppose that $\mathcal{K}$ is nontrapping.
Suppose, also, that $u\in C^\infty$, $u|_{\partial\mathcal{K}}=0$, and
$u(t,x)=0$ for $t<0$.  Then, there is a constant
$C=C_{N_0,\nu_0,\mathcal{K}}$ so that if $u$ vanishes for large $x$
for every fixed $t$,
\begin{multline}\label{3.6}
\sum_{\substack{|\alpha|+\mu\le N_0+\nu_0\\\mu\le \nu_0}} \|\langle
x\rangle^{-(n-1)/4} L^\mu \partial^\alpha u'\|_{L^2_tL^2_x(S_T)} \le
C\int_0^T \sum_{\substack{|\alpha|+\mu\le N_0+\nu_0\\\mu\le \nu_0}}
\|\Box L^\mu \partial^\alpha u(s,\cd)\|_2\:ds
\\+C\sum_{\substack{|\alpha|+\mu\le N_0+\nu_0-1\\\mu\le \nu_0}} \|\Box
L^\mu \partial^\alpha u\|_{L^2_tL^2_x(S_T)}.
\end{multline}
Additionally, 
\begin{multline}\label{3.7}
\sum_{\substack{|\alpha|+\mu\le N_0+\nu_0\\\mu\le \nu_0}} \|\langle
x\rangle^{-(n-1)/4} L^\mu Z^\alpha u'\|_{L^2_tL^2_x(S_T)} \le
C\int_0^T \sum_{\substack{|\alpha|+\mu\le N_0+\nu_0\\\mu\le \nu_0}}
\|\Box L^\mu Z^\alpha u(s,\cd)\|_2\:ds
\\+C\sum_{\substack{|\alpha|+\mu\le N_0+\nu_0-1\\\mu\le \nu_0}} \|\Box
L^\mu Z^\alpha u\|_{L^2_tL^2_x(S_T)}.
\end{multline}
\end{proposition}
While the above proposition is stated for nontrapping geometries, the
same argument will yield estimates for any geometry satisfying the
resolvant bounds used in \cite{B} provided that a sufficient loss of
regularity is allowed for in the right sides.  See \cite{MS}
(Proposition 2.6) for an $n=3$ example.  Additionally, we note that
\eqref{3.5}, \eqref{3.6}, and \eqref{3.7} hold with the weight
$\langle x\rangle^{-(n-1)/4}$ in the left replaced by $\langle
x\rangle^{-1/2-\varepsilon}$ for any $\varepsilon>0$.
We refer the reader to the scaling argument in Keel-Smith-Sogge
\cite{KSS2} (Proposition 2.1) and to the application of such estimates 
in Metcalfe-Sogge-Stewart
\cite{MSS} (Proposition 2.3).  In the sequel, as in \cite{Metcalfe}, we will only
require the estimates as stated.

\noindent{\em Proof of Proposition \ref{prop3.2}:}
Let us prove only \eqref{3.6} as \eqref{3.7} follows from the same
arguments.  

Since the better estimates \eqref{3.local} hold when the norm in the left
is taken over $S_T\cap \{|x|<2\}$, it suffices to consider the norm in the left
over $S_T\cap \{|x|\ge 2\}$.  To do this, we fix $\beta\in C^\infty(\R^n)$ satisfying
$\beta(x)\equiv 1$, $|x|\ge 2$ and $\beta(x)\equiv 0$, $|x|<3/2$.  Since we are assuming
$\mathcal{K}\subset \{|x|<1\}$, it follows that $v=\beta u$, which is equal to $u$ over
$|x|\ge 2$, solves the Minkowski wave equation
$$\Box v = \beta \Box u - 2\nabla\beta\cdot\nabla_x u - (\Delta \beta)u$$
with vanishing initial data.  
%We write $v=v_1+v_2$ where $\Box v_1 = \beta \Box u$,
%$\Box v_2 = -2\nabla\beta\cdot\nabla_x u - (\Delta\beta)u$, and both have zero initial data.  If we 
%apply \eqref{3.5} to $v_1$, we see that this term's contribution 
%is dominated by the right side of \eqref{3.6}. 
%
%It remains to study $v_2$.  Since $\Box v_2$ vanishes unless $3/2<|x|<2$, it follows from \eqref{3.5} 
%that
%$$\sum_{\substack{|\alpha|+\mu\le N_0+\nu_0\\\mu\le \nu_0}} \|\langle x\rangle^{-(n-1)/4} \partial^\alpha
%v_2'\|_{L^2_tL^2_x(S_T)} \le C \sum_{\substack{|\alpha|+\mu\le N_0+\nu_0\\\mu\le\nu_0}} \|L^\mu \partial^\alpha
%u'\|_{L^2_tL^2_x(S_T\cap \{|x|<2\})}.$$
%Here, we have used the fact that the Dirichlet boundary condition allows us to control $u$ locally by $u'$.  The
%bound for the contribution due to $v_2$, thus, follows from \eqref{3.local} which completes the argument.\qed
Here, we apply \eqref{3.5} with $F$ replaced by $\beta\Box u$ and $G$ replaced by $-2\nabla\beta\cdot \nabla_x u
-(\Delta \beta)u$.  It is essential to note that $G$ vanishes unless $|x|<2$.  Thus, by \eqref{3.7}, we have that
\begin{multline*}
\sum_{\substack{|\alpha|+\mu\le N_0+\nu_0\\\mu\le\nu_0}}\|\langle x\rangle^{-(n-1)/4} L^\mu \partial^\alpha v'
\|_{L^2_tL^2_x(S_T)} \le C \int_0^T \sum_{\substack{|\alpha|+\mu\le N_0+\nu_0\\\mu\le \nu_0}} \|L^\mu Z^\alpha
\Box u(s,\cd)\|_2\:ds 
\\+C\sum_{\substack{|\alpha|+\mu\le N_0+\nu_0\\\mu\le\nu_0}} \|L^\mu \partial^\alpha u'\|_{L^2_tL^2_x(S_T\cap
\{|x|<2\})}.
\end{multline*}
Here, we have used the fact that the Dirichlet boundary condition allows us to control $u$ locally by $u'$.
The bound for the last term on the right follows from \eqref{3.local}, which completes the proof.\qed

%%%%%%%%%%%%%%%%%%%%%%%%%%%%%%%%%%%%%%%%%%%%%%%%%%%%%%%%%%%%%%%%
%%%%%%%%%%%%%% SECTION 4 %%%%%%%%%%%%%%%%%%%%%%%%%%%%%%%%%%%%%%%
%%%%%%%%%%%%%%%%%%%%%%%%%%%%%%%%%%%%%%%%%%%%%%%%%%%%%%%%%%%%%%%%

\newsection{Sobolev-type Estimates}

In the sequel, we will require a number of Sobolev-type estimates.  
These are useful for establishing pointwise decay estimates that we be
required in the continuity argument.

We begin with a now standard weighted Sobolev estimate (see \cite{knull}).
\begin{lemma}\label{lemma4.1}
Suppose that $h\in C^\infty(\R^n)$.  Then, for $R\ge 1$,
\begin{equation}\label{4.1}
\|h\|_{L^\infty(R/2<|x|<R)}\le CR^{-(n-1)/2} \sum_{|\alpha|+|\beta|\le
(n+2)/2} \|\Omega^\alpha \partial^\beta_x h\|_{L^2(R/4<|x|<2R)},
\end{equation}
and
\begin{equation}\label{4.2}
\|h\|_{L^\infty(R-1<|x|<R)} \le C R^{-(n-1)} \sum_{|\alpha|+|\beta|\le
n} \|\Omega^\alpha \partial^\beta_x h\|_{L^1(R-2<|x|<R+1)}.
\end{equation}
\end{lemma}

Next, we will need the following estimates for the boundaryless case.
The first is due to Klainerman-Sideris \cite{KS} and says that if
$g\in C^\infty_0(\R_+\times \R^n)$, then
\begin{equation}\label{4.3}
\|\langle t-r\rangle \partial^2 g(t,\cd)\|_2 \le C\sum_{|\alpha|\le 1}
\|\Gamma^\alpha g'(t,\cd)\|_2 + C\|\langle t+r\rangle (\partial_t^2-\Delta)g(t,\cd)\|_2
\end{equation}
where $\Gamma=\{L,Z\}$.
This was shown in \cite{KS} for the $n=3$ case, but the proof is
clearly valid for any $n\ge 2$.  We also have the related estimate
\begin{equation}\label{4.4}
r^{(n/2)-1} \langle t-r\rangle |\partial g(t,x)|\le C\sum_{|\alpha|\le
n/2} \|Z^\alpha \partial g(t,\cd)\|_2 + C\sum_{|\alpha|\le n/2}
\|\langle t-r\rangle Z^\alpha \partial^2 u(t,\cd)\|_2.
\end{equation}
This bound was shown in Hidano \cite{Hidano}.  It is a generalization
of the $n=3$ bound of \cite{HY}.  The latter follows easily from an
estimate of Sideris \cite{Si}.

If we argue as in \cite{MNS} (Lemma 4.2, Lemma 4.3), the above
estimates can be extended to the exterior domain as follows.
\begin{lemma}\label{lemma4.2}
Suppose that $u\in C^\infty_0(\R\times\ext)$ vanishes for
$x\in\partial\mathcal{K}$.  Then, if $|\alpha|=M$ and $\nu$ are fixed,
\begin{multline}\label{4.5}
\|\langle t-r\rangle L^\nu Z^\alpha \partial^2 u(t,\cd)\|_2\le
C\sum_{\substack{|\beta|+\mu\le M+\nu+1\\\mu\le \nu+1}}\|L^\mu Z^\beta
u'(t,\cd)\|_2
\\+C\sum_{\substack{|\beta|+\mu\le M+\nu\\\mu\le \nu}} \|\langle
t+r\rangle L^\mu Z^\beta (\partial_t^2-\Delta) u(t,\cd)\|_2 +
(1+t)\sum_{\mu\le \nu} \|L^\mu u'(t,\cd)\|_{L^2(|x|<2)}
\end{multline}
and
\begin{multline}\label{4.6}
r^{(n/2)-1} \langle t-r\rangle |\partial L^\nu Z^\alpha u(t,x)|\le
C\sum_{\substack{|\beta|+\mu\le M+\nu+(n/2)+1\\\mu\le
\nu+1}} \|L^\mu Z^\beta u'(t,\cd)\|_2
\\+C\sum_{\substack{|\beta|+\mu\le M+\nu+(n/2)\\\mu\le \nu}}
\|\langle t+r\rangle L^\mu Z^\beta (\partial_t^2-\Delta)u(t,\cd)\|_2
+C(1+t)\sum_{\mu\le\nu} \|L^\mu u'(t,\cd)\|_{L^2(|x|<2)}.
\end{multline}
\end{lemma}

%%%%%%%%%%%%%%%%%%%%%%%%%%%%%%%%%%%%%%%%%%%%%%%%%%%%%%%%%%%%%%%%%%%%%%
%%%%%%%%%%%%%%%%%%% SECTION 5%%%%%%%%%%%%%%%%%%%%%%%%%%%%%%%%%%%%%%%%
%%%%%%%%%%%%%%%%%%%%%%%%%%%%%%%%%%%%%%%%%%%%%%%%%%%%%%%%%%%%%%%%%%%%%%
\newsection{Boundary Term Estimates}

In the sequel, we will need to control boundary terms such as those
that appear in \eqref{2.8}, \eqref{4.5}, and \eqref{4.6}.  We will
first need a result for solutions to free wave equations.

\begin{lemma}\label{lemma5.1}
Suppose $n\ge 4$, and suppose that $u_0\in C_0^\infty(\R_+\times\R^n)$
is a solution to the boundaryless wave equation $\Box u_0=G$ with
vanishing initial data.  Then,
\begin{multline}\label{5.1}
\int_0^t \|u(s,\cd)\|_{L^2(|x|<3)}\:ds \le C\int_0^t \|G(s,\cd)\|_2\:ds
%\\+C\int_0^t\int_0^s
%\|G(\tau,\cd)\|_{L^2(||y|-(s-\tau)|<10)}\:d\tau\:ds 
\\+ C\int_0^t \int \sum_{|\alpha|+|\beta|\le n} |\Omega^\alpha
\partial_y^\beta G(s,y)|\:\frac{dy\:ds}{|y|^{(n-1)/2}}.
\end{multline}
\end{lemma}

\noindent{\em Proof of Lemma \ref{lemma5.1}:} Using cutoffs, it
suffices to consider the solution $u(s,\cd)$ in three cases: (1)
$G(\tau,y)$ is supported in $|y|<10$, (2) $G(\tau,y)$ is supported in
$||y|-(s-\tau)|<10$, and (3) $G(\tau,y)$ vanishes unless $|y|>8$ and
$||y|-(s-\tau)|>8$.

The first two cases are handled quite easily.  In the first, we can
use the local energy decay \eqref{1.7} to see that
$$\|u(s,\cd)\|_{L^2(|x|<3)} \le C\int_0^s \frac{1}{(1+s-\tau)^{n-1}}
\|G(\tau,\cd)\|_2\:d\tau.$$
For the second case, we have
$$\|u(s,\cd)\|_{L^2(|x|<3)}\le
C\|u(s,\cd)\|_{L^{\frac{2n}{n-2}}(|x|<3)} \le C\int_0^s
\|G(\tau,\cd)\|_{L^2(||y|-(s-\tau)| <10)}\:d\tau$$
by Sobolev estimates and the energy inequality.

Thus, we only need to establish a bound in the third case.  Here, we
use the fact that
$$u(s,x)=\int_0^s R(s-\tau,\cd)*G(\tau,\cd) \:d\tau$$
where
$$R(t,x)=\lim_{\varepsilon \searrow 0} c_n' \text{Im}
(|x|^2-(t-i\varepsilon)^2)^{-(n-1)/2}.$$
See, e.g., Taylor \cite{Taylor2} p.222.  From this, it clearly follows
that
$$|u(s,x)|\le \int_0^s \int \frac{1}{((s-\tau)^2-|x-y|^2)^{(n-1)/2}}
|G(\tau,y)|\:dy\:d\tau.$$
By support considerations, the right side is bounded by
$$\int_0^s \int \frac{1}{\langle s-\tau-|y|\rangle^{(n-1)/2}}
|G(\tau,y)|\:\frac{dy\:d\tau}{|y|^{(n-1)/2}}$$
when $|x|<3$.

It follows then that
\begin{multline*}
\|u(s,\cd)\|_{L^2(|x|<3)} \le C\int_0^s \frac{1}{(1+s-\tau)^{n-1}}
\|G(\tau,\cd)\|_2\:d\tau
\\+C\int_0^s \|G(\tau,\cd)\|_{L^2(||y|-(s-\tau)|<10)}\:d\tau
\\+ C\int_0^s \int \frac{1}{\langle s-\tau-|y|\rangle^{(n-1)/2}} |G(\tau,y)|\:\frac{dy\:d\tau}{|y|^{(n-1)/2}}.
\end{multline*}
Thus, it is clear that upon integration, we have
\begin{multline}\label{5.2}
\int_0^t \|u(s,\cd)\|_{L^2(|x|<3)}\:ds \le C\int_0^t \|G(s,\cd)\|_2\:ds
\\+C\int_0^t\int_0^s
\|G(\tau,\cd)\|_{L^2(||y|-(s-\tau)|<10)}\:d\tau\:ds + C\int_0^t \int |G(s,y)|\:\frac{dy\:ds}{|y|^{(n-1)/2}}
\end{multline}
provided $n\ge 4$.

The second term on the right side of \eqref{5.2} is dominated by
$$C\int_0^t\int_0^s (s-\tau)^{(n-1)/2}
\|G(\tau,\cd)\|_{L^\infty(||y|-(s-\tau)|<10)}\:d\tau\:ds.$$
By \eqref{4.2}, this is in turn controlled by
$$C\int_0^t \int_0^s \int_{||y|-(s-\tau)|<20}
 \sum_{|\alpha|+|\beta|\le n} |\Omega^\alpha
\partial^\beta G(\tau,y)|\:\frac{dy}{|y|^{(n-1)/2}}\:d\tau\:ds.$$
Since the sets $\{(\tau,y)\,:\,||y|-(j-\tau)|<20\}$, $j=0,1,2,\dots$
have finite overlap, we conclude that this is bounded by
$$C\int_0^t \int \sum_{|\alpha|+|\beta|\le n} |\Omega^\alpha
\partial_x^\beta G(s,y)|\:\frac{dy\:ds}{|y|^{(n-1)/2}}.$$
With this bound, \eqref{5.1} follows immediately from \eqref{5.2}.\qed

We now show how this yields our desired estimates for solutions to
Dirichlet wave equations.  The following is a generalized version of
(2.32) in \cite{MS}.
\begin{lemma}\label{lemma5.2}
Suppose $n\ge 4$, and suppose that $u\in C_0^\infty(\R\times\ext)$
vanishes for $x\in\partial\mathcal{K}$.  Then, if $N_0$ and $\nu\le 1$
are fixed,
\begin{multline}\label{5.3}
\int_0^t \sum_{\substack{|\alpha|+\mu\le
N_0+\nu_0\\\mu\le\nu_0}}\|L^\mu \partial^\alpha
u'(s,\cd)\|_{L^2(|x|<1)}\:ds\le C\int_0^t
\sum_{\substack{|\alpha|+\mu\le N_0+\nu_0+1\\\mu\le \nu_0}} \|L^\mu
\partial^\alpha \Box u(s,\cd)\|_2\:ds
%\\+C\int_0^t\int_0^s \sum_{\substack{|\alpha|+\mu\le
%N_0+\nu_0+1\\\mu\le \nu_0}} \|L^\mu \partial^\alpha\Box
%u(\tau,\cd)\|_{L^2(||y|-(s-\tau)|\le 10)}\:d\tau\:ds
\\+C\int_0^t \int \sum_{\substack{|\alpha|+\mu\le
N_0+\nu_0+n+1\\\mu\le\nu_0}} |L^\mu Z^\alpha \Box
u(s,\cd)|\:\frac{dy\:ds}{|y|^{(n-1)/2}}.
\end{multline}
\end{lemma}

\noindent{\em Proof of Lemma \ref{lemma5.2}:}  Here, we examine two
cases separately: (1) $\Box u(s,y)=0$ for $|y|\ge 4$, and (2) $\Box
u(s,y)=0$ for $|y|\le 3$.  For the former case, we have
\begin{equation}\label{5.4}
\int_0^t \sum_{\substack{|\alpha|+\mu\le N_0+\nu_0\\\mu\le \nu_0}}
\|L^\mu \partial^\alpha u'(t,\cd)\|_{L^2(|x|<1)} \le C\int_0^t
\sum_{\substack{|\alpha|+\mu\le N_0+\nu_0\\\mu\le\nu_0}} \|L^\mu
\partial^\alpha \Box u(s,\cd)\|_2\:ds.
\end{equation}
Indeed, \eqref{1.7} yields
\begin{multline*}
\sum_{\substack{j+\mu\le N_0+\nu_0\\\mu\le \nu_0}} \|\langle
t\rangle^\mu \partial_t^\mu \partial^j_t u'(t,\cd)\|_{L^2(|x|<1)} 
\\\le
C \int \frac{1}{(1+t-s)^{n-1-\nu_0}} \sum_{\substack{|\alpha|+\mu\le
N_0+\nu_0\\\mu\le\nu_0}} \|L^\mu \partial^\alpha
 \Box u(s,\cd)\|_{L^2(|x|<4)} \:ds.
\end{multline*}
Thus, by elliptic regularity (see Lemma 2.3 of \cite{MS}), it follows
that
\begin{multline}\label{5.5}
\sum_{\substack{|\alpha|+\mu\le N_0+\nu_0\\\mu\le\nu_0}} \|L^\mu
\partial^\alpha u'(t,\cd)\|_{L^2(|x|<1)} \le
C\sum_{\substack{|\alpha|+\mu\le N_0+\nu_0-1\\\mu\le\nu_0}} \|L^\mu
\partial^\alpha \Box u(t,\cd)\|_2
\\+C\int_0^t \frac{1}{(1+t-s)^{n-1-\nu_0}}
\sum_{\substack{|\alpha|+\mu\le N_0+\nu_0\\\mu\le\nu_0}} \|L^\mu
\partial^\alpha \Box u(s,\cd)\|_2\:ds.
\end{multline}
This clearly implies \eqref{5.4} for $\nu_0\le 1$ and $n\ge 4$.

In the second case, the case that $\Box u$ vanishes near the obstacle,
we write $u=u_0+u_r$ where $u_0$ solves the boundaryless wave equation
$\Box u_0=\Box u$ with zero initial data.  Fixing $\beta\in
C^\infty(\R^n)$ satisfying $\beta(x)\equiv 1$, $|x|<2$, and
$\beta(x)\equiv 0$ for $|x|>3$, we set $\tilde{u}=\beta u_0+u_r$.
Clearly, $u=\tilde{u}$ for $|x|<2$, and $\tilde{u}$ solves
$$\Box \tilde{u}=-2\nabla\beta\cdot\nabla_x u-(\Delta\beta)u_0$$
which is supported in $|x|<3$.  Thus, from \eqref{5.4}, it follows
that
$$\int_0^t \sum_{\substack{|\alpha|+\mu\le N_0+\nu_0\\\mu\le\nu_0}}
\|L^\mu \partial^\alpha u'(s,\cd)\|_{L^2(|x|<1)} \:ds \le C\int_0^t
\sum_{\substack{|\alpha|+\mu \le N_0+\nu_0 + 1\\\mu\le \nu_0}} \|L^\mu
\partial^\alpha u_0\|_{L^2(|x|<3)}\:ds.$$
Since $[\Box, L]=2\Box$ and $[\Box,\partial]=0$, \eqref{5.3}
is a consequence of \eqref{5.1}. \qed

%%%%%%%%%%%%%%%%%%%%%%%%%%%%%%%%%%%%%%%%%%%%%%%%%%%%%%%%%%%%%%%%
%%%%%%%%%%%%%%%%% SECTION 6 %%%%%%%%%%%%%%%%%%%%%%%%%%%%%%%%%%%%
%%%%%%%%%%%%%%%%%%%%%%%%%%%%%%%%%%%%%%%%%%%%%%%%%%%%%%%%%%%%%%%%
\newsection{Proof of Theorem \ref{theorem1.1}}

In this section, we prove the global existence theorem, Theorem
\ref{theorem1.1}, when $n=4$.  Straightforward modifications will
yield the general case $n\ge 4$.  We take $N=101$ in the smallness
hypothesis \eqref{1.4}; this, however, is not optimal.

The proof of global existence will rely on the following standard
local existence theorem.
\begin{theorem}\label{theorem6.1}
Suppose that $f$ and $g$ are as in Theorem \ref{theorem1.1} with $N\ge
(3n+6)/2$ in \eqref{1.4} if $n$ is even, $N\ge (3n+3)/2$ if $n$ is
odd.  Then there is a $T>0$ so that the initial value problem
\eqref{1.1} with this initial data has a $C^2$ solution satisfying
$$u\in L^\infty([0,T];H^N(\ext))\cap C^{0,1}([0,T];H^{N-1}(\ext)).$$
The supremum of such $T$ is equal to the supremum of all $T$ such that
the initial value problem has a $C^2$ solution with $\partial^\alpha
u$ bounded for $|\alpha|\le 2$.  Also, one can take $T\ge 2$ if
$\|f\|_{H^N}+\|g\|_{H^{N-1}}$ is sufficiently small.
\end{theorem}

This esssentially follows from the local existence results Theorem 9.4
and Lemma 9.6 of Keel-Smith-Sogge \cite{KSS}.  The latter were only
stated for diagonal single-speed systems; however, since the proof
relied only on energy estimates, it extends to the multi-speed,
non-diagonal case if the symmetry assumptions \eqref{1.3} are
satisfied.

Next, in order to avoid dealing with difficulties involving the
compatibility conditions for the Cauchy data, it is convenient to
follow the example of Keel-Smith-Sogge \cite{KSS3} and reduce to an
equivalent equation with vanishing initial data.  We first note that
if the initial data satisfy \eqref{1.4} with $\varepsilon$
sufficiently small, then we can find a solution $u$ to \eqref{1.1} on
a set of the form $0<ct<|x|$ where $c=5\max_I c_I$, and this solution
satisfies
\begin{equation}\label{6.1}
\sup_{0<t<\infty} \sum_{|\alpha|\le 101} \|\langle
x\rangle^{|\alpha|}\partial^\alpha u(t,\cd)\|_{L^2(\ext\,:\,
|x|>ct)}\le C_0\varepsilon.
\end{equation}
Rather than providing unnecessary technicalities, we refer the reader
to \cite{KSS3}, \cite{MS}, or \cite{MNS}.

We will use this local solution $u$ to allow us to restrict to the
case of vanishing Cauchy data.  To do so, we fix a cutoff function
$\chi\in C^\infty(\R)$ satisfying $\chi(s)\equiv 1$ for $s\le
\frac{1}{2c}$ and $\chi(s)\equiv 0$ for $s>\frac{1}{c}$, and we set
$$u_0(t,x)=\eta(t,x)u(t,x),\quad \eta(t,x)=\chi(|x|^{-1}t).$$
Assuming, as we may, that $0\in\mathcal{K}$, we have that $|x|$ is
bounded below on the complement of $\mathcal{K}$ and the function
$\eta(t,x)$ is smooth and homogeneous of degree $0$ in $(t,x)$.  Note
that by \eqref{4.1} and \eqref{6.1}, it follows that there is an
absolute constant $C_1$ so that
\begin{multline}\label{6.2}
(1+t+|x|)\sum_{\mu+|\alpha|\le 98} |L^\mu Z^\alpha u_0(t,x)|
\\+\sum_{\mu+|\alpha|+|\beta|\le 101} \|\langle t+r\rangle^{|\beta|}
L^\mu Z^\alpha \partial^\beta u_0(t,\cd)\|_2\le C_1\varepsilon.
\end{multline}
Notice this also implies that
$$\sum_{\mu+|\alpha|+|\beta|\le 100} \|\langle x\rangle^{-3/4}
L^\mu Z^\alpha u'_0\|_{L^2(S_t)}$$
is $O(\varepsilon)$.  

Since
$$\Box u_0=\eta Q(du,d^2u)+[\Box,\eta]u,$$
$u$ solves $\Box u=Q(du,d^2u)$ for $0<t<T$ if and only if $w=u-u_0$
solves
\begin{equation}\label{6.3}
\begin{cases}
\Box w = (1-\eta)Q(du,d^2u)-[\Box,\eta]u\\
w|_{\partial\mathcal{K}}=0\\
w(t,x)=0,\quad t\le 0
\end{cases}
\end{equation}
for $0<t<T$.

If we let $v$ be the solution of the linear equation
\begin{equation}\label{6.4}
\begin{cases}
\Box v = -[\Box,\eta]u\\
v|_{\partial\mathcal{K}}=0\\
v(t,x)=0,\quad t\le 0,
\end{cases}
\end{equation}
then we will show that \eqref{6.1} implies that there is another
absolute constant $C_2$ so that
\begin{equation}\label{6.5}
\sum_{\mu+|\alpha|\le 99} \|L^\mu Z^\alpha v'(t,\cd)\|_2 +
\sum_{\mu+|\alpha|\le 100} \|\langle x\rangle^{-3/4} L^\mu Z^\alpha
v'(t,\cd)\|_{L^2_tL^2_x(S_t)}\le C_2\varepsilon.
\end{equation}

Indeed, we can examine the first term in \eqref{6.5} using the
standard energy integral method.  Doing so, we see that
\begin{multline*}
\partial_t \sum_{\mu+|\alpha|\le 99} \|L^\mu Z^\alpha v'(t,\cd)\|_2^2 
\\\le C\Bigl(\sum_{\mu+|\alpha|\le 99}\|L^\mu Z^\alpha
v'(t,\cd)\|_2\Bigr)\Bigl(\sum_{\mu+|\alpha|\le 99} \|L^\mu Z^\alpha
\Box v(s,\cd)\|_2\Bigr)
\\+C\sum_{\mu+|\alpha|\le 99}
\Bigl|\int_{\partial\mathcal{K}}\partial_0 L^\mu Z^\alpha
v(t,\cd)\nabla L^\mu Z^\alpha v(t,\cd)\cdot n\:d\sigma\Bigr|,
\end{multline*}
where $n$ is the outward normal at a given point on
$\partial\mathcal{K}$.  Since $\mathcal{K}\subset\{|x|<1\}$ and since
$\Box v= -[\Box,\eta]u$, it follows that the right side of the
equation above is dominated by
\begin{multline*}
C\Bigl(\sum_{\mu+|\alpha|\le 99} \|L^\mu Z^\alpha
v'(t,\cd)\|_2\Bigr)\Bigl(\sum_{\mu+|\alpha|\le 99} \|L^\mu Z^\alpha
[\Box,\eta]u(s,\cd)\|_2\Bigr)
\\+C\int_{\{x\in\extiv\,:\,|x|<1\}} \sum_{\mu+|\alpha|\le 100} |L^\mu
Z^\alpha v'(t,x)|^2\:dx.
\end{multline*}
Integrating in $t$, this yields
\begin{multline}\label{6.6}
\sum_{\mu+|\alpha|\le 99}\|L^\mu Z^\alpha v'(t,\cd)\|_2^2 \le
C\Bigl(\int_0^t \sum_{\mu+|\alpha|\le 99} \|L^\mu Z^\alpha
[\Box,\eta]u(s,\cd)\|_2\:ds\Bigr)^2\\+C\int_0^t \sum_{\mu+|\alpha|\le
100} \|L^\mu Z^\alpha v'(s,\cd)\|^2_{L^2(|x|<1)}\:ds.
\end{multline}
The last term in \eqref{6.6} is dominated by the square of the second
term in the left side of \eqref{6.5}.

Thus, by \eqref{3.7}, it follows that the square of the left side of
\eqref{6.5} is controlled by
$$C\Bigl(\int_0^t \sum_{\mu+|\alpha|\le 100} \|L^\mu Z^\alpha
[\Box,\eta]u(s,\cd)\|_2\:ds\Bigr)^2 + C\sum_{\mu+|\alpha|\le 99}
\|L^\mu Z^\alpha [\Box,\eta]u\|^2_{L^2_tL^2_x(S_t)}.$$
Both of these terms are $O(\varepsilon^2)$ by \eqref{6.1}, which
completes the proof of \eqref{6.5}.

Using this, we are now ready to set up the continuity argument which
will complete the proof of Theorem \ref{theorem1.1}.  If
$\varepsilon>0$ is as above, we shall assume that we have a $C^2$
solution of \eqref{6.3} for $0\le t\le T$ satisfying the lossless
estimates
\begin{align}
\sum_{\substack{|\alpha|+\mu\le 52\\\mu\le 1}} \|L^\mu Z^\alpha
w'(t,\cd)\|_2 + \sum_{\substack{|\alpha|+\mu \le 53\\\mu\le 1}}
\|\langle x\rangle^{-3/4} L^\mu Z^\alpha w'\|_{L^2_tL^2_x(S_t)}&\le
A_0\varepsilon \label{6.7}\\
(1+t+r)\sum_{|\alpha|\le 40} |Z^\alpha w'(t,x)|&\le B_1\varepsilon,\label{6.8}
\end{align}
as well as, the lossy higher order estimates
\begin{align}
\sum_{|\alpha|\le 100} \|\partial^\alpha w'(t,\cd)\|_2&\le
B_2\varepsilon (1+t)^{1/40}\label{6.9}\\
\sum_{\substack{|\alpha|+\mu\le 81\\\mu\le 2}} \|L^\mu Z^\alpha
w'(t,\cd)\|_2 &\le B_3\varepsilon (1+t)^{1/20}\label{6.10}\\
\sum_{\substack{|\alpha|+\mu\le 79\\\mu\le 2}} \|\langle x\rangle^{-3/4}
L^\mu Z^\alpha w'\|_{L^2_tL^2_x(S_t)} 
&\le B_4\varepsilon (1+t)^{1/20}.\label{6.11}
\end{align}
Here, as before, the $L^2_x$-norms are taken over $\extiv$, and the
weighted $L^2_tL^2_x$ norms are over $S_t=[0,t]\times\extiv$.

In \eqref{6.7}, we may take $A_0=10C_2$ where $C_2$
is the constant in \eqref{6.5}.
Clearly, if $\varepsilon$ is sufficiently small, then all of these
estimates hold for $T=2$ by Theorem \ref{theorem6.1}.  With this in
mind, we shall then prove that, for $\varepsilon>0$ smaller than some
number depending on $B_1,\dots,B_4$,
\begin{enumerate}
\item[$(i.)$] \eqref{6.7} is valid with $A_0$ replaced by $A_0/2$,
\item[$(ii.)$] \eqref{6.8}-\eqref{6.11} are consequences of \eqref{6.7}.
\end{enumerate}
It will then follow from the local existence theorem that a solution
exists for all $t>0$ if $\varepsilon$ is small enough.

\noindent{\em Proof of $(i.)$:}  

Since $v$ satisfies the better bound \eqref{6.5}, it suffices to show 
\begin{equation}\label{6.12}
\sum_{\substack{|\alpha|+\mu\le 52\\\mu\le 1}} \|L^\mu Z^\alpha
(w-v)'(t,\cd)\|^2_2 + \sum_{\substack{|\alpha|+\mu\le 53\\\mu\le 1}}
\|\langle x\rangle^{-3/4} L^\mu Z^\alpha (w-v)'\|^2_{L^2_tL^2_x(S_t)}\le C\varepsilon^4.
\end{equation}
Using the energy integral method as in the proof of \eqref{6.5}, it
follows that the first term on the left side of \eqref{6.12} is
controlled by
$$
C\Bigl(\int_0^t \sum_{\substack{|\alpha|+\mu\le 52\\\mu\le 1}} \|L^\mu
Z^\alpha \Box u(s,\cd)\|_2\:ds\Bigr)^2 + C
\sum_{\substack{|\alpha|+\mu\le 53\\\mu\le 1}} \|L^\mu Z^\alpha
(w-v)'\|^2_{L^2_tL^2_x(S_t \cap\{|x|<1\})}
$$
since $\Box (w-v)=(1-\eta)\Box u$.  Thus, by \eqref{3.7}, the left
side of \eqref{6.12} is controlled by
\begin{equation}\label{6.13}
C\Bigl(\int_0^t \sum_{\substack{|\alpha|+\mu\le 53\\\mu\le 1}} \|L^\mu
Z^\alpha \Box u(s,\cd)\|_2\:ds\Bigr)^2 +
C\sum_{\substack{|\alpha|+\mu\le 52\\\mu\le 1}} \|L^\mu Z^\alpha \Box
u\|^2_{L^2_tL^2_x(S_t)}.
\end{equation}
We will show that the first term in \eqref{6.13} is
$O(\varepsilon^4)$.  The same techniques can be applied to get the
bound for the second term.

We begin by noting that for $|\beta|+\nu\le 53$, $\nu\le 1$, we have
\begin{multline}\label{6.14}
|L^\nu Z^\beta \Box u(s,y)| \le
|u'(s,y)|\sum_{\substack{|\alpha|+\mu\le 53\\\mu\le 1}}|L^\mu Z^\alpha
\partial^2 u(s,y)| 
\\+
\sum_{\substack{|\alpha|+\mu\le
30\\\mu\le 1}} |L^\mu Z^\alpha u'(s,y)|
\sum_{\substack{|\alpha|+\mu\le 53\\\mu\le 1}} |L^\mu Z^\alpha u'(s,y)|.
\end{multline}
Thus, by \eqref{4.1} and \eqref{4.5}, we have
\begin{multline}\label{6.15}
\sum_{\substack{|\alpha|+\mu\le 53\\\mu\le 1}} \|L^\mu Z^\alpha \Box
u(s,\cd)\|_2 
\\\le \frac{C}{(1+s)^{3/2}} \sum_{|\alpha|\le 3} \|Z^\alpha
u'(s,\cd)\|_2 \sum_{\substack{|\alpha|+\mu\le 54\\\mu\le 1}} \|L^\mu
Z^\alpha u'(s,\cd)\|_{L^2(|y|>\tilde{c}s/2)}
\\+\frac{C}{1+s}\sum_{|\alpha|\le 3} \|\langle y\rangle^{-3/2} Z^\alpha
u'(s,\cd)\|_{L^2(|y|<\tilde{c}s/2)}
\sum_{\substack{|\alpha|+\mu\le 54\\\mu\le 2}} \|L^\mu Z^\alpha
u'(s,\cd)\|_2
\\+\frac{C}{1+s}\sum_{|\alpha|\le 3} \|\langle y\rangle^{-3/2}
Z^\alpha u'(s,\cd)\|_{L^2(|y|<\tilde{c}s/2)}
\sum_{\substack{|\alpha|+\mu\le 53\\\mu\le 1}} \|\langle s+|y|\rangle
L^\mu Z^\alpha \Box u(s,\cd)\|_2 
\\+C \sum_{|\alpha|\le 3} \|\langle y\rangle^{-3/2} Z^\alpha u'(s,\cd)\|_2 
\sum_{\mu\le 1} \|L^\mu u'(s,\cd)\|_{L^2(|y|<1)}
\\+C\sum_{\substack{|\alpha|+\mu\le 53\\\mu\le 1}} \|\langle
y\rangle^{-3/4} L^\mu Z^\alpha u'(s,\cd)\|^2_2.
\end{multline}
Here $\tilde{c}=(1/2)\min_I c_I$.
The right side of \eqref{6.15} is in turn bounded by
\begin{multline*}
\le \frac{C}{(1+s)^{3/2}} \sum_{\substack{|\alpha|+\mu\le 54\\\mu\le
2}}\|L^\mu Z^\alpha u'(s,\cd)\|_2^2  +C \sum_{\substack{|\alpha|+\mu\le
53\\\mu\le 1}} \|\langle y\rangle^{-3/4} L^\mu Z^\alpha
u'(s,\cd)\|_2^2 
\\+\frac{C}{(1+s)^2} \sum_{\substack{|\alpha|+\mu\le 53\\\mu\le 1}}
\|L^\mu Z^\alpha u'(s,\cd)\|_2^4
\\+C \sum_{|\alpha|\le 3} \|Z^\alpha u'(s,\cd)\|_2
\sum_{\substack{|\alpha|+\mu\le 53\\\mu\le 1}} \|L^\mu Z^\alpha \Box
u(s,\cd)\|_{L^2(|y|<\tilde{c}s/2)}
\end{multline*}
using \eqref{4.1} in the third term on the right side of \eqref{6.15}
if $|y|>\tilde{c}s/2$.
By \eqref{6.2} and \eqref{6.7}, the last term can be absorbed into the left side of
\eqref{6.15} if $\varepsilon$ is small enough.  Thus, we see that
\begin{multline}\label{6.16}
\int_0^t \sum_{\substack{|\alpha|+\mu\le 53\\\mu\le 1}} \|L^\mu
Z^\alpha \Box u(s,\cd)\|_2\:ds \le C\int_0^t \frac{1}{(1+s)^{3/2}}
\sum_{\substack{|\alpha|+\mu\le 54\\\mu\le 2}} \|L^\mu Z^\alpha
u'(s,\cd)\|_2^2\:ds
\\+C\sum_{\substack{|\alpha|+\mu\le 53\\\mu\le 1}} \|\langle
y\rangle^{-3/4} L^\mu Z^\alpha u'\|^2_{L^2_tL^2_x(S_t)}
+C\int_0^t \frac{1}{(1+s)^2} \sum_{\substack{|\alpha|+\mu\le
53\\\mu\le 1}} \|L^\mu Z^\alpha u'(s,\cd)\|_2^4 \:ds.
\end{multline}
The first and third terms of \eqref{6.16} are easily seen to be
$O(\varepsilon^2)$ by \eqref{6.2} and \eqref{6.10}.  The second term is also
$O(\varepsilon^2)$ by \eqref{6.2} and \eqref{6.7}.  This shows that the
first term in \eqref{6.13} is $O(\varepsilon^4)$ as desired.  Since
same argument can be used to establish that the second term in \eqref{6.13} is also
$O(\varepsilon^4)$, we have \eqref{6.12} which completes the proof of $(i.)$.

%%%%%%%%%%%%%%%%%%%%%%%%%%%%%%%%%%%%%%%%%%%%%%
\noindent{\em Proof of $(ii.)$:}

In this section, we complete the proof of Theorem \ref{theorem1.1} by
showing that \eqref{6.8}-\eqref{6.11} are consequences of
\eqref{6.7}.  Our first task will be to establish \eqref{6.8}.  Over
$|x|\ge \tilde{c}t/2$, the bound follows from \eqref{4.1} and
\eqref{6.7}.  Over $|x|<\tilde{c}t/2$, we apply \eqref{4.6} to bound
the left side of \eqref{6.8} by
\begin{multline}\label{6.17}
C\sum_{\substack{|\alpha|+\mu\le 43\\\mu\le 1}} \|L^\mu Z^\alpha
w'(t,\cd)\|_2
+C\sum_{|\alpha|\le 42} \|\langle t+r\rangle Z^\alpha \Box w(t,\cd)\|_2 
\\+C(1+t)\|w'(t,\cd)\|_{L^2(|x|<1)}.
\end{multline}
The first term of \eqref{6.17} is $O(\varepsilon)$ by \eqref{6.7}.
When $\Box w$ in the second term is replaced by $\Box u_0=\eta
Q(du,d^2u) +[\Box,\eta]u$, which is supported in $|x|>ct$,
this term is seen to be $O(\varepsilon)$ by a Sobolev estimate and 
\eqref{6.1}.  When $\Box w$ is replaced by $\Box u$ in the second
term, we see that it is bounded by
\begin{multline*}
C \sum_{|\alpha|\le 43} \|Z^\alpha u'(t,\cd)\|_2\:
\sup_x\Bigl((1+t+|x|) \sum_{|\alpha|\le 40} |Z^\alpha u'(t,\cd)|\Bigr)
\\\le C\varepsilon^2 + C\varepsilon \sup_x \Bigl((1+t+|x|)
\sum_{|\alpha|\le 40} |Z^\alpha w'(t,\cd)|\Bigr)
\end{multline*}
by \eqref{6.2} and \eqref{6.7}.  For $\varepsilon$ sufficiently small,
the second of these terms can be absorbed into the left side of \eqref{6.8}.
Since $\Box w = \Box u - \Box
u_0$, this establishes the desired control on the second term in
\eqref{6.17}.  

It remains to bound the last term in \eqref{6.17}.  By the fundamental
theorem of calculus, we see that it is dominated by
$$C\int_0^t \sum_{|\alpha|+\mu\le 1} \|L^\mu \partial^\alpha
w'(s,\cd)\|_{L^2(|x|<1)}\:ds.$$
By \eqref{5.3}, this is in turn controlled by
\begin{equation}\label{6.18}
C\int_0^t \sum_{\substack{|\alpha|+\mu\le 2\\\mu\le 1}} \|L^\mu
\partial^\alpha \Box w(s,\cd)\|_2\:ds
%\\+C\int_0^t \int_0^s \sum_{\substack{|\alpha|+\mu\le 2\\\mu\le 1}}
%\|L^\mu \partial^\alpha \Box
%w(\tau,\cd)\|_{L^2(||y|-(s-\tau)|<10)}\:d\tau\:ds
+C\int_0^t \int \sum_{\substack{|\alpha|+\mu\le 6\\\mu\le 1}} |L^\mu
Z^\alpha \Box w(s,y)|\:\frac{dy\:ds}{|y|^{3/2}}.
\end{equation}
When $\Box w$ is replaced by $\Box u_0$, it follows from \eqref{6.1}
that both of these terms are $O(\varepsilon)$.  In the remaining case,
when $\Box w$ is replaced by $\Box u$, it follows from \eqref{4.1}
that both terms in \eqref{6.18} are controlled by
$$C\sum_{\substack{|\alpha|+\mu\le 6\\\mu\le 1}} \|\langle
x\rangle^{-3/4} L^\mu Z^\alpha u'\|^2_{L^2_tL^2_x(S_t)}.$$
This completes the proof of \eqref{6.8} since the above is
$O(\varepsilon)$ by \eqref{6.2} and \eqref{6.7}.

With \eqref{6.8} established, the remainder of the proof follows very
similarly to that in \cite{MS}.  The main exception is how we deal
with the boundary term in \eqref{2.8}.  We will only provide a sketch
of the arguments that follow exactly as in \cite{MS}.  The reader may
also wish to refer to \cite{MNS}.

Let us begin with the proof of \eqref{6.9}.  In the notation of \S 2,
we have $(\Box_\gamma u)^I=\sum_{\substack{0\le j,k\le
4\\1\le J,K\le D}} A^{I,jk}_{JK}\partial_j u^J \partial_k u^K$ and
$\gamma^{IJ,jk}=-\sum_{\substack{0\le l\le 4\\1\le K\le D}}
B^{IJ,jk}_{K,l} \partial_l u^K$.  Notice that by \eqref{6.2} and \eqref{6.8}
\begin{equation}\label{6.19}
\|\gamma'(s,\cd)\|_{\infty}\le\frac{C\varepsilon}{(1+s)}.
\end{equation}

In order to prove \eqref{6.9}, we first estimate the energy of
$\partial_t^j u$ for $j\le M\le 100$ using induction on $M$.  By
\eqref{2.5} and \eqref{6.19}, we have
\begin{equation}\label{6.20}
\partial_t E^{1/2}_M(u)(t)\le C\sum_{j\le M} \|\Box_\gamma
\partial^j_t u(t,\cd)\|_2 + \frac{C\varepsilon}{1+t} E^{1/2}_M(u)(t).
\end{equation}
Since it follows from elliptic regularity and \eqref{6.8} that
\begin{multline*}
\sum_{j\le M}\|\Box_\gamma \partial^j_t u(t,\cd)\|_2\le
\frac{C\varepsilon}{1+t}\sum_{j\le M} \|\partial^j_t u'(t,\cd)\|_2
\\+C\sum_{|\alpha|\le M-41,|\beta|\le M/2} \|\partial^\alpha
u'(t,\cd)\partial^\beta u'(t,\cd)\|_2,
\end{multline*}
we obtain
\begin{equation}\label{6.21}
\partial_t E_M^{1/2}(u)(t)\le \frac{C\varepsilon}{1+t}E^{1/2}_M(u)(t)
+ C\sum_{|\alpha|\le M-41, |\beta|\le M/2} \|\partial^\alpha
u'(t,\cd)\partial^\beta u'(t,\cd)\|_2
\end{equation}
since $E_M^{1/2}(u)(t)\approx \sum_{j\le M} \|\partial_t^j
u'(t,\cd)\|_2$ for $\varepsilon$ small.

When $M=40$, the last term in \eqref{6.21} drops out.  Thus, since
\eqref{1.4} implies that $E_{100}^{1/2}(u)(0) \le C\varepsilon$, Gronwall's
inequality yields
\begin{equation}\label{6.22}
\sum_{j\le 40} \|\partial_t^j u'(t,\cd)\|_2\le C\varepsilon
(1+t)^{C\varepsilon}.\end{equation}
For $M>40$, we have to deal with the last term in \eqref{6.21}.  By
\eqref{4.1}, this term is bounded by
$$C\sum_{|\alpha|\le \max(M-38,3+M/2)} \|\langle x\rangle^{-3/4}Z^\alpha
u'(t,\cd)\|^2_2.$$
Thus, \eqref{6.21} and Gronwall's inequality yield,
\begin{equation}\label{6.23}
E_M^{1/2}(u)(t)\le C(1+t)^{C\varepsilon} \Bigl[\varepsilon +
\sum_{|\alpha|\le \max(M-38,3+M/2)} \|\langle x\rangle^{-3/4} Z^\alpha
u'\|^2_{L^2_tL^2_x(S_t)}\Bigr].
\end{equation}
If we use \eqref{6.22} and \eqref{6.23}, 
\begin{equation}\label{6.24}
E_{100}^{1/2}(u)(t)\le C\varepsilon (1+t)^{C\varepsilon+\sigma}
\end{equation}
would follow for arbitrarily small $\sigma>0$ from a simple induction
argument and the following lemma.  At every step of the induction, we
are using the fact that bounds on $E_M^{1/2}(u)(t)$ yield bounds on
$\sum_{|\alpha|\le M} \|\partial^\alpha u'(t,\cd)\|_2$ by elliptic
regularity.

\begin{lemma}\label{lemma6.2}
Under the above assumptions, if $M\le 100$ and
\begin{multline}\label{6.25}
\sum_{|\alpha|\le M} \|\partial^\alpha u'(t,\cd)\|_2 +
\sum_{|\alpha|\le M-3} \|\langle x\rangle^{-3/4} \partial^\alpha
u'\|_{L^2_tL^2_x(S_t)} +\sum_{|\alpha|\le M-4} \|Z^\alpha
u'(t,\cd)\|_2
\\+\sum_{|\alpha|\le M-6} \|\langle x\rangle^{-3/4} Z^\alpha
u'\|_{L^2_tL^2_x(S_t)} \le C\varepsilon (1+t)^{C\varepsilon+\sigma}
\end{multline}
with $\sigma>0$, then there is a constant $C'$ so that
\begin{multline}\label{6.26}
\sum_{|\alpha|\le M-2} \|\langle x\rangle^{-3/4} \partial^\alpha
u'\|_{L^2_tL^2_x(S_t)} + \sum_{|\alpha|\le M-3} \|Z^\alpha
u'(t,\cd)\|_2 
\\+\sum_{|\alpha|\le M-5} \|\langle x\rangle^{-3/4} Z^\alpha
u'\|_{L^2_tL^2_x(S_t)} \le C'\varepsilon (1+t)^{C'\varepsilon+C'\sigma}.
\end{multline}
\end{lemma}

\noindent{\em Proof of Lemma \ref{lemma6.2}:} We start by estimating
the first term in the left side of \eqref{6.26}.  By \eqref{6.2}, \eqref{6.5},
\eqref{3.6}, and the fact that $\Box (w-v)=(1-\eta)\Box u$, this is dominated by
$$C\varepsilon+C\sum_{|\alpha|\le M-2} \int_0^t \|\partial^\alpha \Box
u(s,\cd)\|_2\:ds + C\sum_{|\alpha|\le M-3} \|\partial^\alpha \Box
u\|_{L^2_tL^2_x(S_t)}.$$
If $M\le 40$, we can use \eqref{6.2}, \eqref{6.8}, and \eqref{6.25} to see that the
last two terms are $\le C\varepsilon^2 (1+t)^{C\varepsilon+\sigma}$.
If $40<M\le 100$, we can repeat the proof of \eqref{6.23} to conclude
that the are
\begin{align*}
&\le C\varepsilon^2 (1+t)^{C\varepsilon+\sigma} + C\sum_{|\alpha|\le
\max(M-40,3+M/2)} \|\langle x \rangle^{-3/4} Z^\alpha
u'\|^2_{L^2_tL^2_x(S_t)}
\\&\quad\quad\quad +C\sup_{0\le s\le  t}\Bigl(\sum_{|\alpha|\le M-6}
\|Z^\alpha u'(s,\cd)\|_2\Bigr) \sum_{|\alpha|\le \max(M-40,3+M/2)}
\|\langle x\rangle^{-3/4} Z^\alpha u'\|_{L^2_tL^2_x(S_t)}\\
&\le C\varepsilon^2 (1+t)^{2C\varepsilon + 2\sigma},
\end{align*}
using the inductive hypothesis \eqref{6.25} and the fact that
$\max(M-40,3+M/2)\le M-6$ if $M\ge 40$.

We next establish the bound for the second term in the left of \eqref{6.26}
using \eqref{2.10}.  With $Y_{M-3,0}(t)$ as in \eqref{2.9}, it
follows as in the proof of \eqref{6.25} that
$$\sum_{|\alpha|\le M-3} \|\Box_\gamma Z^\alpha u(t,\cd)\|_2 \le
\frac{C\varepsilon}{1+t} Y^{1/2}_{M-3,0}(t)+C\sum_{|\beta|\le M-40}
\|\langle x\rangle^{-3/4} Z^\beta u'(t,\cd)\|_2^2$$
using \eqref{4.1}, \eqref{6.2}, and \eqref{6.8}.  Plugging this into \eqref{2.10}, we
have
\begin{multline*}
\partial_t Y_{M-3,0}(t)\le
\frac{C\varepsilon}{1+t}Y_{M-3,0}(t)+C\sum_{|\beta|\le M-40} \|\langle
x\rangle^{-3/4} Z^\beta u'(t,\cd)\|_2^2
\\+C\sum_{|\alpha|\le M-2} \|\langle x\rangle^{-3/4}\partial^\alpha
u'(t,\cd)\|_2^2.
\end{multline*}
By Gronwall's inequality and the fact that $\sum_{|\alpha|\le M-3}
\|Z^\alpha u'(t,\cd)\|_2^2\le CY_{M-3,0}(t)$ for $\varepsilon$ small
enough, this yields
\begin{multline*}
\sum_{|\alpha|\le M-3} \|Z^\alpha u'(t,\cd)\|_2^2 \le
C(1+t)^{C\varepsilon} \Bigl(\varepsilon^2 + C\sum_{|\beta|\le M-40}
\|\langle x\rangle^{-3/4} Z^\beta u'\|^2_{L^2_tL^2_x(S_t)}
\\+C\sum_{|\alpha|\le M-2} \|\langle x\rangle^{-3/4} \partial^\alpha
u'\|^2_{L^2_tL^2_x(S_t)}\Bigr).
\end{multline*}
The last term is bounded by the right side of \eqref{6.26} using
the previous step.  For the second term in the right, we can apply the
inductive hypothesis \eqref{6.25} which yields \eqref{6.26}.

Using \eqref{3.7}, this in turn implies that the third term in the
left of \eqref{6.26} satisfies the bounds, which completes the
proof. \qed

This proves \eqref{6.24}.  By elliptic regularity and \eqref{6.2},
\eqref{6.9} follows.  It also follows from the lemma that
\begin{multline}\label{6.27}
\sum_{|\alpha|\le 98} \|\langle x\rangle^{-3/4} \partial^\alpha
w'\|_{L^2_tL^2_x(S_t)} + \sum_{|\alpha|\le 97} \|Z^\alpha
w'(t,\cd)\|_2 \\+\sum_{|\alpha|\le 95} \|\langle x\rangle^{-3/4}
Z^\alpha w'\|_{L^2_tL^2_x(S_t)} \le C\varepsilon
(1+t)^{C\varepsilon+\sigma}.
\end{multline}
Here and in what follows $\sigma$ denotes a small constant that must
be taken to be larger at each occurence.

We now proceed to the proof of the estimates involving powers of $L$.
We first estimate $L^\nu \partial^\alpha u'$ in $L^2$.  We then
obtain \eqref{6.11} and \eqref{6.12} for this $\nu$ using an inductive
argument similar to Lemma \ref{lemma6.2}.

The main part of the next step is to show that
\begin{equation}\label{6.28}
\sum_{\substack{|\alpha|+\mu\le 92\\\mu\le 1}} \|L^\mu \partial^\alpha
u'(t,\cd)\|_2\le C\varepsilon (1+t)^{C\varepsilon+\sigma}.
\end{equation}
For this, we shall use \eqref{2.8}.  We must first establish
\eqref{2.7} for $N_0+\nu_0\le 92$, $\nu_0=1$.  Arguing as in the proof
of \eqref{6.23}, which uses \eqref{4.1}, \eqref{6.8}, and elliptic
regularity, we get that for $M\le 92$
\begin{multline*}
\sum_{\substack{j+\mu\le M\\\mu\le 1}} \Bigl(\|\tilde{L}^\mu
\partial^j_t \Box_\gamma u(t,\cd)\|_2 + \|[\tilde{L}^\mu \partial^j_t,
\Box-\Box_\gamma] u(t,\cd)\|_2\Bigr)
\\\le \frac{C\varepsilon}{1+t}\sum_{\substack{j+\mu\le M\\\mu\le 1}}
\|\tilde{L}^\mu \partial_t^j u'(t,\cd)\|_2
\\+C\sum_{|\alpha|\le M-41} \|\langle x\rangle^{-3/4} L \partial^\alpha
u'(t,\cd)\|_2 \sum_{|\alpha|\le 95} \|\langle x\rangle^{-3/4} Z^\alpha
u'(t,\cd)\|_2
\\+C\sum_{|\alpha|\le \max(M,3+M/2)} \|\langle x\rangle^{-3/4} Z^\alpha
u'(t,\cd)\|_2^2.
\end{multline*}
Thus, we obtain \eqref{2.7} with $\delta=C\varepsilon$ and
$$H_{1,M-1}(t)=C\sum_{|\alpha|\le M-41} \|\langle x\rangle^{-3/4}
L\partial^\alpha u'(t,\cd)\|_2^2 + C\sum_{|\alpha|\le 95} \|\langle
x\rangle^{-3/4} Z^\alpha u'(t,\cd)\|^2_2.$$
Since \eqref{1.4} gives $\Bigl(\int e_0(\tilde{L}^\mu\partial_t^j
u)(0,x)\:dx\Bigr)^{1/2} \le C\varepsilon$ if $\mu+j\le
100$, it follows from \eqref{2.8} and \eqref{6.27} that for $M\le 92$
\begin{multline}\label{6.29}
\sum_{\substack{|\alpha|+\mu\le M\\\mu\le 1}} \|L^\mu \partial^\alpha
u'(t,\cd)\|_2 \le C\varepsilon (1+t)^{C\varepsilon+\sigma} +
C(1+t)^{C\varepsilon} \sum_{|\alpha|\le M-41} \|\langle x\rangle^{-3/4}
L\partial^\alpha u'\|^2_{L^2_tL^2_x(S_t)}
\\+C(1+t)^{C\varepsilon} \int_0^t \sum_{|\alpha|\le M+1}
\|\partial^\alpha u'(s,\cd)\|_{L^2(|x|<1)}\:ds.
\end{multline}
By \eqref{6.2} and \eqref{5.3}, this last integral is dominated by
$\varepsilon \log(2+t)$ plus
\begin{multline}\label{6.30}
\int_0^t \sum_{|\alpha|\le M+1} \|\partial^\alpha
w'(s,\cd)\|_{L^2(|x|<1)} \:ds 
\le C\int_0^t \sum_{|\alpha|\le M+2} \|\partial^\alpha \Box
w(s,\cd)\|_2\:ds
\\+C\int_0^t \int \sum_{|\alpha|\le M+6} |\partial^\alpha \Box w(s,\cd)|\:\frac{dy\:ds}{|y|^{3/2}}.
\end{multline}
When $w$ is replaced by $u_0$, both of these terms are
$O(\varepsilon)$ by \eqref{6.1}.  Since $\Box w= \Box u -\Box u_0$, it
suffices to consider the case that $w$ is replaced by $u$.  In this
case, the right side of \eqref{6.30} is controlled by
$$C\sum_{|\alpha|\le 95} \|\langle x\rangle^{-3/4} Z^\alpha
u'\|^2_{L^2_tL^2_x(S_t)} + C\sum_{|\alpha|\le 98} \|\langle
x\rangle^{-3/4} \partial^\alpha u'\|_{L^2_tL^2_x(S_t)}.$$
Both of these terms are in turn $\le C \varepsilon
(1+t)^{C\varepsilon+\sigma}$ by \eqref{6.2} and \eqref{6.27}.

Therefore, by \eqref{6.29}, we have that
\begin{multline*}
\sum_{\substack{|\alpha|+\mu\le M\\\mu\le 1}} \|L^\mu \partial^\alpha
u'(t,\cd)\|_2 \le C\varepsilon (1+t)^{C\varepsilon + \sigma}
\\+C(1+t)^{C\varepsilon} \sum_{|\alpha|\le M-41} \|\langle
x\rangle^{-3/4} L \partial^\alpha u'\|^2_{L^2(S_t)}.
\end{multline*}
This gives the desired bound when $M\le 40$.  Since the analog of
Lemma \ref{lemma6.2} is valid when $M=100$ is replaced by $M=92$ and
$u$ is replaced by $Lu$, we get \eqref{6.28} by a simple induction
argument.  This same induction also yields, as in the case of no
$L$'s,
\begin{multline}\label{6.31}
\sum_{\substack{|\alpha|+\mu\le 90\\\mu\le 1}} \|\langle
x\rangle^{-3/4} L^\mu \partial^\alpha w'\|_{L^2_tL^2_x(S_t)} +
\sum_{\substack{|\alpha|+\mu\le 89\\\mu\le 1}} \|L^\mu Z^\alpha
w'(t,\cd)\|_2 
\\+\sum_{\substack{|\alpha|+\mu\le 87\\\mu\le 1}} \|\langle
x\rangle^{-3/4} L^\mu Z^\alpha w'\|_{L^2_tL^2_x(S_t)} \le C\varepsilon
(1+t)^{C\varepsilon + C\sigma}.
\end{multline}
Repeating this argument for $L^2 Z^\alpha u'$, it in turn follows from
\eqref{6.28} and \eqref{6.31} that
$$
\sum_{\substack{|\alpha|+\mu\le 81\\\mu\le 2}} \|L^\mu
Z^\alpha w'(t,\cd)\|_2 
\\+\sum_{\substack{|\alpha|+\mu\le 79\\\mu\le 2}} \|\langle
x\rangle^{-3/4} L^\mu Z^\alpha w'\|_{L^2_tL^2_x(S_t)}\le C\varepsilon
(1+t)^{C\varepsilon+C\sigma},
$$
which implies \eqref{6.10} and \eqref{6.11}.  This completes the proof
of $(ii.)$, and hence the proof of Theorem \ref{theorem1.1}.

%%%%%%%%%%%%%%%%%%%%%%%%%%%%%%%%%%%%%%%%%%%%%%%%%%%%%%%%%%%%%%%%%
%%%%%%%%%%%%%%%%%%%%%% BIBLIOGRAPHY %%%%%%%%%%%%%%%%%%%%%%%%%%%%%
%%%%%%%%%%%%%%%%%%%%%%%%%%%%%%%%%%%%%%%%%%%%%%%%%%%%%%%%%%%%%%%%%


\begin{thebibliography}{MA}
%\bibitem{AY} R. Agemi and K. Yokoyama: {\em The null condition and
%global existence of solutions to systems of wave equations with
%different speeds}, Advances in Nonlinear Partial Differential
%Equations and Stochastics, (1998), 43--86.
%\bibitem{burq} N. Burq: {\em D\'ecroissance de l\'energie locale
%de l'\'equation des ondes pour le probl\`me ext\'erieur et absence
%de r\'sonance au voisinage du r\'eel}, Acta Math. {\bf 180}
%(1998), 1--29.
\bibitem{B} N. Burq: {\em Global Strichartz estimates for nontrapping
geometries: A remark about an artical by H. Smith and C. Sogge}, to
appear, Comm. Partial Differentail Equations.
%\bibitem{christ} D. Christodoulou: {\em Global solutions of
%nonlinear hyperbolic equations for small initial data}, Comm. Pure
%Appl. Math. {\bf 39}, (1986), 267-282.
%\bibitem{D} P. S. Datti:  {\em Nonlinear wave equations in exterior domains},
% Nonlinear Anal. {\bf 15} (1990), 321--331.
% \bibitem{gilbarg} D. Gilbarg and N. Trudinger:
%{\em Elliptic partial differential equations of second order},
%Springer, Second Ed., Third Printing, 1998.
\bibitem{Hayashi} N. Hayashi: {\em Global existence of small solutions
to quadratic nonlinear wave equations in an exterior domain},
J. Funct. Anal., {\bf 131} (1995), 302--344.
\bibitem{Hidano} K. Hidano: {\em An elementary proof of global or
almost global existence for quasi-linear wave equations},
preprint.
\bibitem{HY} K. Hidano and K. Yokoyama: {\em A remark on the almost
global existence theorems of Keel, Smith, and Sogge}, preprint.
\bibitem{HY2} K. Hidano and K. Yokoyama: {\em A new proof of the
global existence theorem of Klainerman}, preprint.
\bibitem{H} L. H\"ormander:
{\em Lectures on nonlinear hyperbolic equations}, Springer-Verlag,
Berlin, 1997.
%\bibitem{h2}
%L. H\"ormander: {\em $L^1, L^\infty$ estimates for the wave
%operator}, Analyse Mathematique et Applications, Gauthier-Villars,
%Paris, 1988, pp. 211-234.
\bibitem{Ikawa1} M. Ikawa: {\em Decay of solutions of the wave
equation in the exterior of two convex bodies}, Osaka J. Math.
{\bf 19} (1982), 459--509.
\bibitem{Ikawa2} M. Ikawa: {\em Decay of solutions of the wave
equation in the exterior of several convex bodies}, Ann. Inst. Fourier
(Grenoble), {\bf 38} (1988), 113--146.
\bibitem{KSS} M. Keel, H.
Smith, and C. D. Sogge: {\em Global existence for a
 quasilinear wave equation
outside of star-shaped domains}, J. Funct. Anal. {\bf{189}},
(2002), 155--226.
\bibitem{KSS2} M. Keel, H. Smith, and C. D. Sogge: {\em Almost
global existence for some semilinear wave equations}, J.
D'Analyse, {\bf 87} (2002), 265--279.
\bibitem{KSS3}M. Keel, H.
Smith, and C. D. Sogge: {\em Almost global existence for 
quasilinear wave equations in three space dimensions},
J. Amer. Math. Soc. {\bf 17} (2004), 109--153.
%\bibitem{K} S. Klainerman:
%{\em Uniform decay estimates and the Lorentz invariance
% of the classical wave equation},
%Comm. Pure Appl. Math. {\bf 38}(1985), 321--332.
\bibitem{knull} S. Klainerman:
{\em The null condition and global existence to nonlinear wave
equations}, Lectures in Applied Math. \textbf{23} (1986),
293--326.
\bibitem{KS} S. Klainerman and T. Sideris: {\em On almost global existence for
nonrelativistic wave equations in 3d} Comm. Pure Appl. Math. {\bf
49}, (1996), 307--321.
%\bibitem{KY} K. Kubota and K. Yokoyama: {\em Global existence of
%classical solutions to systems of nonlinear wave equations with
%different speeds of propagation}, Japan. J. Math. {\bf 27} (2001), 113--202.
\bibitem{LMP} P. D. Lax, C. S. Morawetz, and R. S. Phillips:
{\em Exponential decay of solutions of the  wave equation in the
exterior of a star-shaped obstacle}, Comm. Pure Appl. Math. {\bf
16} (1963), 477--486.
\bibitem{LP} P. D. Lax and R. S. Phillips:
{\em Scattering theory, revised edition}, Academic Press, San
Diego, 1989.
\bibitem{Mel} R. B. Melrose: {\em Singularities and energy decay
of acoustical scattering}, Duke Math. J. {\bf 46} (1979), 43--59.
\bibitem{Metcalfe} J. Metcalfe: {\em Global existence for semilinear
wave equations exterior to nontrapping obstacles}, Houston
J. Math. {\bf 30} (2004), 259--281.
\bibitem{Metcalfe2} J. Metcalfe: {\em Global Strichartz estimates for
solutions of the wave equation exterior to convex obstacles}, Trans. Amer.
Math. Soc. {\bf 356} (2004), 4839--4855.
\bibitem{MS} J. Metcalfe and C. D. Sogge: {\em Hyperbolic trapped
rays and global existence of quasilinear wave equations}, Invent. Math. {\bf 159}
(2005), 75--117.
\bibitem{MNS} J. Metcalfe, M. Nakamura, and C. D. Sogge: {\em Global
existence of solutions to multiple speed systems of quasilinear wave
equations in exterior domains}, Forum Math. {\bf 17} (2005), 133--168.
\bibitem{MSS} J. Metcalfe, C. D. Sogge, and A. Stewart: {\em Nonlinear
hyperbolic equations in infinite homogeneous waveguides}, Comm. Partial
Differential Equations, to appear.
\bibitem{M} C. S. Morawetz: {\em The decay of solutions of the exterior
initial-boundary problem for the wave equation}, Comm. Pure Appl.
Math. {\bf 14} (1961), 561--568.
\bibitem{MRS} C. S. Morawetz, J. Ralston, and W. Strauss:
{\em Decay of solutions of the wave equation outside nontrapping
obstacles}, Comm. Pure Appl. Math. {\bf 30} (1977), 447--508.
\bibitem{ralston} J. V. Ralston: {\em Solutions of the wave
equation with localized energy}, Comm. Pure Appl. Math., {\bf 22)}
(1969), 807--923.
\bibitem{Ralston2} J. Ralston: {\em Note on the decay of acoustic
waves}, Duke Math. J., {\bf 46} (1979), 799--804.
\bibitem{tsutsumi} Y. Shibata and Y. Tsutsumi:
{\em On a global existence theorem of small amplitude solutions
for nonlinear wave equations in an exterior domain}, Math. Z. {\bf
191} (1986), 165-199.
\bibitem{Si} T. Sideris: {\em Nonresonance and global existence of prestressed
nonlinear elastic waves},
 Ann. of Math. {\bf 151} (2000), 849--874.
\bibitem{Si2} T. Sideris: {\em The null condition and global existence of
nonlinear elastic waves.},
Inven. Math. {\bf 123} (1996), 323--342.
\bibitem{Si3} T. Sideris, S.Y. Tu: {\em Global existence
for systems of nonlinear wave equations in 3D with multiple
speeds.}, SIAM J. Math. Anal. {\bf 33} (2001), 477--488.
\bibitem{SS} H. Smith and C. D. Sogge: {\em Global Strichartz estimates for
nontrapping perturbations of the Laplacian}, Comm. Partial
Differential Equations {\bf 25} (2000), 2171--2183.
\bibitem{S}  C. D. Sogge: {\em Lectures on nonlinear wave equations},
International Press, Cambridge, MA, 1995.
\bibitem{So2} C.D. Sogge: {\em Global existence for nonlinear
wave equations with multiple speeds}, Harmonic analysis at Mount
Holyoke (South Hadley, MA, 2001), 353--366, Contemp. Math., 320,
Amer. Math. Soc., Providence, RI, 2003.
\bibitem{Strauss} W. Strauss: {\em Dispersal of waves vanishing on the
boundary of an exterior domain}, Comm. Pure Appl. Math., {\bf 28}
(1975), 265--278.
\bibitem{Taylor} M. Taylor: {\em Grazing rays and reflection of
singularities of solutions to wave equations}, Comm. Pure Appl. Math.,
{\bf 29} (1976), 1--38.
\bibitem{Taylor2} M. Taylor: {\em Partial Differential Equations I},
Springer-Verlag, New York, NY, 1996.
\bibitem{Vainberg} B. R. Vainberg: {\em On the short wave asymptotic
behavior of solutions of stationary problems and the asymptotic
behavior as $t\to\infty$ of solutions of non-stationary problems},
Russian Math Surveys, {\bf 30:2} (1975), 1--55.
\bibitem{Y} K. Yokoyama: {\em Global existence of classical solutions to
systems of wave equations with critical nonlinearity in three
space dimensions} J. Math. Soc. Japan {\bf 52} (2000), 609--632.
\end{thebibliography}
\end{document}